\documentclass[12pt,a4paper,twoside]{article}

\newcommand{\pageformat}[6]{\setlength{\hoffset}{-1in}
                  \setlength{\voffset}{-1in}
                  \addtolength{\hoffset}{#5}
                            \addtolength{\voffset}{#6}
                            \setlength{\oddsidemargin}{#1}
                            \setlength{\evensidemargin}{#2}
                            \setlength{\textwidth}{\paperwidth}
                  \addtolength{\textwidth}{-\oddsidemargin}
                  \addtolength{\textwidth}{-\evensidemargin}
                  \addtolength{\textwidth}{-\marginparsep}
                  \addtolength{\textwidth}{-\marginparwidth}
                            \setlength{\topmargin}{#3}
                            \setlength{\textheight}{\paperheight}
                  \addtolength{\textheight}{-\topmargin}
                  \addtolength{\textheight}{-\headheight}
                  \addtolength{\textheight}{-\headsep}
                  \addtolength{\textheight}{-\footskip}
                  \addtolength{\textheight}{-#4}}
\pageformat{2cm}{3cm}{25mm}{20mm}{1pt}{0pt}

\usepackage{ifthen}
\newboolean{article}
    \setboolean{article}{true}
\newboolean{report}
\newboolean{book}
\newboolean{letter}
\newboolean{german}
\newboolean{italian}
\newboolean{nobaselinestretch}
\newboolean{nosectionappendix}
\newboolean{oldtoc}
\newboolean{nosectionequn}
\newboolean{notheorem}

\ifthenelse{\boolean{german}}{
    \usepackage{german}}{}

\usepackage[latin1]{inputenc}

\usepackage{amsmath}
\usepackage{amssymb}
\usepackage[mathscr]{eucal}

\ifthenelse{\boolean{notheorem}}{}{
    \usepackage{theorem}}

\ifthenelse{\boolean{nobaselinestretch}}{}{
    \renewcommand{\baselinestretch}{1.25}}

\newenvironment{env}[2]{\begin{#1}#2\end{#1}}{}
    \newcommand{\beq}[1]{\begin{env}{equation}{#1}}
    \newcommand{\beqn}[1]{\begin{env}{equation*}{#1}}
    \newcommand{\bal}[1]{\begin{env}{align}{#1}}
    \newcommand{\baln}[1]{\begin{env}{align*}{#1}}
    \newcommand{\bga}[1]{\begin{env}{gather}{#1}}
    \newcommand{\bgan}[1]{\begin{env}{gather*}{#1}}
    \newcommand{\bflal}[1]{\begin{env}{flalign}{#1}}
    \newcommand{\bflaln}[1]{\begin{env}{flalign*}{#1}}
    \newcommand{\bmu}[1]{\begin{env}{multline}{#1}}
    \newcommand{\bmun}[1]{\begin{env}{multline*}{#1}}
    \newcommand{\bsp}[1]{\begin{env}{split}{#1}}

    \newcommand{\eeq}{\end{env}}
    \newcommand{\eeqn}{\end{env}}
    \newcommand{\eal}{\end{env}}
    \newcommand{\ealn}{\end{env}}
    \newcommand{\ega}{\end{env}}
    \newcommand{\egan}{\end{env}}
    \newcommand{\eflal}{\end{env}}
    \newcommand{\eflaln}{\end{env}}
    \newcommand{\emu}{\end{env}}
    \newcommand{\emun}{\end{env}}
    \newcommand{\esp}{\end{env}}

\newcommand{\lf}{\vspace{2ex}}

\renewcommand{\bf}[1]{\textbf{#1}}
\renewcommand{\it}[1]{\textit{#1}}

\renewcommand{\sf}[1]{\textsf{#1}}

\renewcommand{\tt}[1]{\texttt{#1}}
\newcommand{\hl}[1]{\bf{\it{#1}}}

\newcommand{\mbf}[1]{\mathbf{#1}}
\newcommand{\msf}[1]{\text{\small$\sf{#1}$}}

\newcommand{\cmc}[1]{\mathcal{#1}}
\newcommand{\eus}[1]{\mathscr{#1}}

\newcommand{\bb}[1]{\mathbb{#1}}

\newcommand{\mscriptsize}[1]{{\setlength{\arraycolsep}{.3ex}\text{\scriptsize$#1$}}}
\newcommand{\mtiny}[1]{{\setlength{\arraycolsep}{.3ex}\text{\tiny$#1$}}}
\newcommand{\nbd}[1]{$#1$\nobreakdash--}
\newcommand{\ol}[1]{\overline{#1}}
\newcommand{\ul}[1]{\underline{#1}}

\newcommand{\ve}{\varepsilon}
\newcommand{\vt}{\vartheta}

\newcommand{\vp}{\varphi}

\newcommand{\Om}{\Omega}

\newcommand{\norm}[1]{\left\lVert#1\right\rVert}

\newcommand{\bnorm}[1]{\bigl\lVert#1\bigr\rVert}

\newcommand{\Bnorm}[1]{\Bigl\lVert#1\Bigr\rVert}

\newcommand{\bfam}[1]{\bigl(#1\bigr)}
\newcommand{\Bfam}[1]{\Bigl(#1\Bigr)}
\newcommand{\AB}[1]{\langle#1\rangle}

\newcommand{\CB}[1]{\{#1\}}
\newcommand{\bCB}[1]{\bigl\{#1\bigr\}}
\newcommand{\BCB}[1]{\Bigl\{#1\Bigr\}}
\newcommand{\SB}[1]{[#1]}

\newcommand{\Matrix}[1]{\begin{pmatrix}#1\end{pmatrix}}

\newcommand{\sMatrix}[1]{\mscriptsize{\Matrix{#1}}}
\newcommand{\tMatrix}[1]{\mtiny{\Matrix{#1}}}

\newcommand{\set}[2][]{
    \ifthenelse{\equal{#1}{}}{
        \CB{#2}}{
        \CB{#1~|~#2}}}
\newcommand{\bset}[2][]{
    \ifthenelse{\equal{#1}{}}{
        \bCB{#2}}{
        \bCB{#1~|~#2}}}
\newcommand{\Bset}[2][]{
    \ifthenelse{\equal{#1}{}}{
        \BCB{#2}}{
        \BCB{#1~\big|~#2}}}
\newcommand{\zero}{\CB{0}}

\DeclareMathOperator{\ls}{\normalfont\msf{span}}
\DeclareMathOperator{\cls}{\ol{\ls}}

\DeclareMathOperator{\id}{\normalfont\msf{id}}

\renewcommand{\ker}{\operatorname{\msf{ker}}}

\newcommand{\C}{\bb{C}}

\newcommand{\E}{\bb{E}}

\newcommand{\R}{\bb{R}}

\newcommand{\cA}{\cmc{A}}
\newcommand{\cB}{\cmc{B}}
\newcommand{\cC}{\cmc{C}}
\newcommand{\cD}{\cmc{D}}

\newcommand{\sB}{\eus{B}}
\newcommand{\sC}{\eus{C}}

\newcommand{\sK}{\eus{K}}

\newcommand{\U}{\mbf{1}}

\ifthenelse{\boolean{nosectionequn}}{}{
    \numberwithin{equation}{section}
    }

\ifthenelse{\boolean{article}\or\boolean{letter}\or\boolean{nosectionequn}}{
    \setboolean{nosectionappendix}{true}}{}
\ifthenelse{\boolean{nosectionappendix}}{}{
    \renewcommand{\appendix}{
        \chapter*{\appendixname}
        \addcontentsline{toc}{chapter}{\appendixname}
        \renewcommand{\thesection}{\Alph{section}}
        \setcounter{section}{0}}}
   
\ifthenelse{\boolean{report}\or\boolean{book}}{
    }{}

\ifthenelse{\boolean{notheorem}}{}{
        \newcommand{\mnname}{Mathematical note.}
        \newcommand{\enname}{End of the note.}
        \newcommand{\definame}{Definition.}
        \newcommand{\propname}{Proposition.}
        \newcommand{\lemname}{Lemma.}
        \newcommand{\exname}{Example.}
        \newcommand{\exername}{Exercise.}
        \newcommand{\remname}{Remark.}
        \newcommand{\obname}{Observation.}
        \newcommand{\thmname}{Theorem.}
        \newcommand{\corname}{Corollary.}
        \newcommand{\proofname}{Proof.}
    \ifthenelse{\boolean{german}}{
        \renewcommand{\mnname}{Mathematische Notiz.}
        \renewcommand{\enname}{Ende der Notiz.}
        \renewcommand{\exname}{Beispiel.}
        \renewcommand{\exername}{Übung.}
        \renewcommand{\remname}{Bemerkung.}
        \renewcommand{\obname}{Beobachtung.}
        \renewcommand{\thmname}{Satz.}
        \renewcommand{\corname}{Korollar.}
        \renewcommand{\proofname}{Beweis.}}{}
    \ifthenelse{\boolean{italian}}{
        \renewcommand{\mnname}{Nota matematica.}
        \renewcommand{\enname}{Fina della nota.}
        \renewcommand{\definame}{Definizione.}
        \renewcommand{\propname}{Proposizione.}
        \renewcommand{\exname}{Esempio.}
        \renewcommand{\exername}{Esercizio.}
        \renewcommand{\remname}{Nota.}
        \renewcommand{\obname}{Osservazione.}
        \renewcommand{\thmname}{Teorema.}
        \renewcommand{\corname}{Corollario.}
        \renewcommand{\proofname}{Dimostrazione.}

       \renewcommand{\appendixname}{Appendice}

       }{}
    \theoremheaderfont{\normalfont\bfseries}
    \theoremstyle{change}
        \theorembodyfont{\rmfamily}
            \newtheorem{emp}{}[section]
                \newcommand{\bemp}[1][]{
                    \begin{emp}\hskip-\labelsep\bf{#1}\hskip\labelsep}
                \newcommand{\eemp}{\end{emp}}
\newtheorem{itemp}[emp]{}
                \newcommand{\bitemp}[1][]{
                    \begin{itemp}\hskip-\labelsep\bf{#1}\hskip\labelsep\normalfont\itshape}
                \newcommand{\eitemp}{\end{itemp}}
            \newtheorem{mn}[emp]{\mnname}
                \newcommand{\bnm}{\begin{mn}~\begin{quotation}\renewcommand{\baselinestretch}{1}\small\noindent\ignorespaces}
                \newcommand{\enm}{\end{quotation}\hfill\bf{\enname}\end{mn}}
            \newtheorem{ex}[emp]{\exname}
                \newcommand{\bex}{\begin{ex}}
                \newcommand{\eex}{\end{ex}}
            \newtheorem{exer}[emp]{\exername}
                \newcommand{\bexer}{\begin{exer}}
                \newcommand{\eexer}{\end{exer}}
            \newtheorem{defi}[emp]{\definame}
                \newcommand{\bdefi}{\begin{defi}}
                \newcommand{\edefi}{\end{defi}}
            \newtheorem{rem}[emp]{\remname}
                \newcommand{\brem}{\begin{rem}}
                \newcommand{\erem}{\end{rem}}
            \newtheorem{ob}[emp]{\obname}
                \newcommand{\bob}{\begin{ob}}
                \newcommand{\eob}{\end{ob}}
        \theorembodyfont{\normalfont\itshape}
            \newtheorem{thm}[emp]{\thmname}
                \newcommand{\bthm}{\begin{thm}}
                \newcommand{\ethm}{\end{thm}}
            \newtheorem{prop}[emp]{\propname}
                \newcommand{\bprop}{\begin{prop}}
                \newcommand{\eprop}{\end{prop}}
            \newtheorem{cor}[emp]{\corname}
                \newcommand{\bcor}{\begin{cor}}
                \newcommand{\ecor}{\end{cor}}
            \newtheorem{lem}[emp]{\lemname}
                \newcommand{\blem}{\begin{lem}}
                \newcommand{\elem}{\end{lem}}
\newenvironment{empn}[1]{\lf\noindent\bf{#1}\ignorespaces\hskip\labelsep}{\lf}
		\newcommand{\bempn}[1]{\begin{empn}{#1}}
		\newcommand{\eempn}{\end{empn}}
		\newcommand{\bitempn}[1]{\begin{empn}{#1}\normalfont\itshape}
		\newcommand{\eitempn}{\end{empn}}
                \newcommand{\bnmn}{\begin{empn}{\mnname}~\begin{quotation}\renewcommand{\baselinestretch}{1}\small\noindent\ignorespaces}
                \newcommand{\enmn}{\end{quotation}\hfill\bf{\enname}\end{empn}}
		\newcommand{\bexn}{\begin{empn}{\exname}}
		\newcommand{\eexn}{\end{empn}}
		\newcommand{\bexern}{\begin{empn}{\exername}}
		\newcommand{\eexern}{\end{empn}}
		\newcommand{\bdefin}{\begin{empn}{\definame}}
		\newcommand{\edefin}{\end{empn}}
		\newcommand{\bremn}{\begin{empn}{\remname}}
		\newcommand{\eremn}{\end{empn}}
		\newcommand{\bobn}{\begin{empn}{\obname}}
		\newcommand{\eobn}{\end{empn}}

\newcommand{\qedsymbol}{~\rule[-0.35mm]{2mm}{2mm}}
    \newcounter{proof}[emp]
    \newenvironment{Proof}[1]{
        \vspace{1ex}
        \renewcommand{\item}[1][\stepcounter{proof}(\roman{proof})]%
            {##1\hskip\labelsep}
        \noindent\textsc{#1\hskip\labelsep}}{
        \nolinebreak\qedsymbol}
    \newcommand{\proof}[1][\proofname]{
        \begin{Proof}{#1}\ignorespaces}
    \newcommand{\qed}{\end{Proof}}
    \newcommand{\noqed}{
        \renewcommand{\qedsymbol}{}
        \end{Proof}}}
    \ifthenelse{\boolean{italian}}{
        \renewcommand{\proofname}{Dimostrazione.}}{}

\usepackage[varg]{txfonts}

\begin{document}

\title{Generalized Unitaries and the Picard Group}
\author{}
\author{
~\\
Michael Skeide\thanks{This work is supported by research funds of the University of Molise.}\\\\
{\small\itshape Dipartimento S.E.G.e S.}\\
{\small\itshape Università degli Studi del Molise}\\
{\small\itshape Via de Sanctis}\\
{\small\itshape 86100 Campobasso, Italy}\\
{\small{\itshape E-mail: \tt{skeide@math.tu-cottbus.de}}}\\
{\small{\itshape Homepage: \tt{http://www.math.tu-cottbus.de/INSTITUT/lswas/\_skeide.html}}}\\
\\
}
\date{November 2005, revised April 2006}

{
\renewcommand{\baselinestretch}{1}
\maketitle



\begin{abstract}
\noindent
After discussing some basic facts about generalized module maps, we use the representation theory of the algebra $\sB^a(E)$  of adjointable operators on a Hilbert \nbd{\cB}module $E$ to show that the quotient of the group of generalized unitaries on $E$ and its normal subgroup of unitaries on $E$ is a subgroup of the group of automorphisms of the range ideal $\cB_E$ of $E$ in $\cB$. We determine the kernel of the canonical mapping into the Picard group of $\cB_E$ in terms of the group of quasi inner automorphisms of $\cB_E$. As a by-product we identify the group of bistrict automorphisms of the algebra of adjointable operators on $E$ modulo inner automorphisms as a subgroup of the (opposite of the) Picard group.
\end{abstract}

}

\vspace{6ex}


\noindent
A generalized unitary on a Hilbert \nbd{\cB}module $E$ is a surjection $u$ on $E$ satisfying $\AB{ux,uy}=\vp(\AB{x,y})$ for some  automorphism $\vp$ of $\cB$. By conjugation with $u$ we define a bistrict automorphism of the algebra $\sB^a(E)$ of all adjointable maps on $E$. By an application of the theory of strict representations we show that the group of bistrict automorphisms modulo the normal subgroup of inner autormorphisms of $\sB^a(E)$ is a subgroup of $Pic(\cB_E)^{op}$, the opposite of the Picard group of the \hl{range ideal} $\sB_E:=\cls\AB{E,E}$, that is, the isomorphism classes of Morita equivalences from $\cB_E$ to $\cB_E$ with the tensor product as group operation. We determine the kernel of the canonical map from the bistrict automorphisms induced by generalized unitaries into the Picard group in terms of the group of generalized unitaries modulo the normal subgroup of unitaries, which turns out to be a subgroup of the group of automorphisms of $\cB_E$.

In Section \ref{genmodmapSEC} we define generalized module maps and analyze their basic properties. In particular, we prove that a \nbd{\vp}linear map factors into the canonical map from $E$ into its extension via $\vp$ and a usual module map. We suppose that much of Section \ref{genmodmapSEC} will be \it{folclore} (except, possibly, the mentioned factorization). In particular, generalized unitaries have been discussed in Bakic and Guljas \cite{BaGu02}. We emphasize, however, that the terminology used in \cite{BaGu02} is different. (What we call a generalized unitary they call just unitary. But, a unitary on a Hilbert module in most other papers is a surjection in $\sB^a(E)$ that preserves inner products. We definitely prefer to follow the usual terminology and not \cite{BaGu02}.) In Section \ref{PicSEC} we use the representation theory of $\sB^a(E)$ to analyze the role of the Picard group for the group of bistrict automorphisms. As the representation theory (Muhly, Skeide and Solel \cite{MSS06}) is rather new, we expect that Section \ref{PicSEC} consists largely of new material. In Section \ref{combSEC}, finally, we put together the results to explain the relation of the group of generalized unitaries and the Picard group.

Our motivation for these notes is to study one-parameter groups of generalized unitaries and to prepair the terrain for the discussion of cocycles of generalized unitaries on $E$ with respect to \nbd{E_0}semigroups on $\sB^a(E)$. While the inclusion of these discussions is completely out of the range of these notes (and far from having reached an end), we found it convenient to present also rather known things in Section \ref{genmodmapSEC} in order to underline the role played by the tensor product. The latter seems to have been neglected so far in literature. Starting from Section \ref{PicSEC} the tensor product becomes indispensable and in the future applications that we have in mind the tensor product will be the key ingredient of the approach. We comment on these ideas in the ends of Section \ref{genmodmapSEC} and of Section \ref{combSEC}.

\section{Generalized module maps}\label{genmodmapSEC}

Let $\vp\colon\cB\rightarrow\cC$ be a homomorphism between \nbd{C^*}algebras $\cB$ and $\cC$. Then the Hilbert \nbd{\cC}mod\-ule $\cC$ (with inner product $\AB{c,c'}=c^*c'$) inherits a left action of elements $b\in\cB$ by setting $bc:=\vp(b)c$. We say $\vp$ is \hl{nondegenerate}, if this left action is nondegenerate, that is, if $\cB\cC$ is total in $\cC$. In this case $\cC$ is a \hl{correspondence} from $\cB$ to $\cC$, that is, a Hilbert \nbd{\cC}module with a nondegenerate representation of $\cB$. We denote that correspondence by $_\vp\cC$.

\bob\label{auob}
If $\cB$ is unital, then nondegeneracy simply means that $\vp$ is unital. If $\cB$ is nonunital, then nondegeneracy is equivalent to saying that the image of any bounded approximate unit $\bfam{u_\lambda}_{\lambda\in\Lambda}$ for $\cB$ converges \nbd{*}strongly in $\sB^a(_\vp\cC)$ to $\id_{_\vp\cC}$ or, equivalently, $\bfam{\vp(u_\lambda)}_{\lambda\in\Lambda}$ is a bounded approximate unit for $\cC$. (This follows by three epsilons from the inequality
\beqn{
\norm{c-u_\lambda c}
~\le~
\norm{c-c_0}+\norm{c_0-u_\lambda c_0}+\norm{u_\lambda(c_0-c)},
}\eeqn
where $c\in\cC$ is arbitrary and $c_0$ is in $\ls\cB\cC$ sufficiently close to $c$, and the observation that $u_\lambda c_0\to c_0$.) In particular, if $\cC$ is unital, then $\vp(u_\lambda)$ converges to $\U_\cC$ in norm, and if, in this case, $\vp$ is injective, then necessarily $\cB$ is also unital.
\eob

Suppose $E$ is a Hilbert \nbd{\cB}module. Then the tensor product $E\odot{_\vp\cC}$ of $E$ and $_\vp\cC$ is a Hilbert \nbd{\cC}module, the \hl{extension} of $E$ by $\cC$ via $\vp$. We observe that $i_\vp\colon xb\mapsto x\odot\vp(b)$ (in particular, $x\mapsto x\odot\U_\cC$ in the unital case) well-defines a map, the \hl{canonical map}, $E\rightarrow E\odot{_\vp\cC}$. Indeed, we verify that $\AB{i_\vp(xb),i_\vp(x'b')}=\vp(b^*)\vp(\AB{x,x'})\vp(b')=\vp(\AB{xb,x'b'})$ so that $i_\vp$ is a contraction on the linear hull of $E\cB$ that extends to a unique contraction defined on all of $E$. When $\cB$ is unital, then we observe more simply that $i_\vp\colon x\mapsto x\odot\U_\cC$ is well-defined and fulfills $i_\vp(xb)=x\odot\vp(b)=i_\vp(x)\vp(b)$.

\bob\label{insurob}
$i_\vp$ is injective, if and only if the restriction of $\vp$ to $\cB_E$ is injective. (To see necessity choose an element $b\in\ker\vp\cap\cB_E$ and an element $x\in E$ such that $xb\ne0$.) Certainly, $i_\vp$ is surjective, if $\vp$ is, but this condition is not necessary. A necessary and sufficient condition is that $\vp(\cB_E)$ be a (right) ideal. (To see necessity suppose that $\vp(\cB_E)$ is not an ideal and, therefore, by \nbd{*}invariance not even a right ideal. This means there are $c\in\cC,b\in\cB_E$ such that $\vp(b)c\notin\vp(\cB_E)$. Therefore, there exists $x\in E$ such that $x\odot c\notin i_\vp(E)$.)
\eob

The two properties of $i_\vp$ mentioned before Observation \ref{insurob} motivate the following definition.

\bdefi
Let $\vp\colon\cB\rightarrow\cC$ be a nondegenerate homomorphism of \nbd{C^*}algebras. Let $E$ and $F$ denote a Hilbert \nbd{\cB}module and a Hilbert \nbd{\cC}module, respectively.
\begin{enumerate}
\item
An additive map $a\colon E\rightarrow F$ is \hl{\nbd{\vp}linear}, if $a(xb)=(ax)\vp(b)$ for all $x\in E,b\in\cB$. More generally, $a$ is a \hl{generalized module map}, if there exists a nondegenerate homomorphism $\vp$ such that $a$ is \nbd{\vp}linear. We denote the spaces of bounded \nbd{\vp}linear maps and of bounded generalized module maps from $E$ to $F$ by $\sB^\vp(E,F)$ and by $\sB^{gen}(E,F)$, respectively.

\item
A map $a\colon E\rightarrow F$ is \hl{\nbd{\vp}adjointable}, if there exists a linear map $a^*\colon F\rightarrow E$ fulfilling $\AB{ax,y}=\vp(\AB{x,a^*y})$ for all $x\in E,y\in F$. More generally, $a$ is \hl{generalized adjointable}, if there exists a nondegenerate homomorphism $\vp$ such that $a$ is \nbd{\vp}adjointable. We denote the spaces of \nbd{\vp}adjointable \nbd{\vp}linear maps and of generalized adjointable generalized module maps from $E$ to $F$ by $\sB^{\vp,a}(E,F)$ and by $\sB^{gen,a}(E,F)$, respectively.

\item
A map $v\colon E\rightarrow F$ is a \hl{\nbd{\vp}isometry}, if $\AB{vx,vy}=\vp(\AB{x,y})$ for all $x,y\in E$. More generally, $v$ is a \hl{generalized isometry}, if there exists a nondegenerate homomorphism $\vp$ such that $v$ is a \nbd{\vp}isometry.

\item
Suppose $\vp$ is an isomorphism. A map $u\colon E\rightarrow F$ is a \hl{\nbd{\vp}unitary}, if it is a surjective \nbd{\vp}isometry. More generally, $u$ is a \hl{generalized unitary}, if there exists an isomorphism $\vp$ such that $u$ is a \nbd{\vp}unitary. We denote by $U^{gen}(E,F)$ the set of generalized unitaries $E\rightarrow F$. In particular, we denote by $U^{gen}(E):=U^{gen}(E,E)$ the \hl{generalized unitary group} of $E$.
\end{enumerate}
Of course, we will use abbreviations like $U^{gen}(E)=U^{gen}(E,E)$ for the generalized unitaries also for all other spaces of maps.
\edefi

Recall that in a Hilbert \nbd{\cB}module we have $\lim_\lambda xu_\lambda=x$ for every approximate unit $\bfam{u_\lambda}_{\lambda\in\Lambda}$ for $\cB$. Taking into account also Observation \ref{auob} we find that a \nbd{\vp}linear map is, in particular, linear.

In the following observations we collect a couple of basic properties. They illustrate that the most useful cases occur when $\vp$ is injective. They also illustrate the useful method to check equality of elements in a (pre-)Hilbert module by comparing their inner products.

\bob\label{isolinob}
Every \nbd{\vp}isometry is \nbd{\vp}linear. This follows, because in the pre-Hilbert $\vp(\cB)$--mod\-ule $vE$ by
\beqn{
\AB{vx,v(y+zb)}
~=~
\vp(\AB{x,y+zb})
~=~
\vp(\AB{x,y})+\vp(\AB{x,z})\vp(b)
~=~
\AB{vx,vy+(vz)\vp(b)}
}\eeqn
all the inner products of the elements $v(y+zb)$ and $vy+(vz)\vp(b)$ with $vx$ $(x\in E)$ coincide. Actually, $vE$ is complete. To see this it is sufficient to find a homomorphism $\vt_v\colon\sB^a(E)\rightarrow\sB^a(vE)$ such that
\beqn{
\Phi_v
\colon
\Matrix{b&y^*\\x&a}
~\longrightarrow~
\Matrix{\vp(b)&(vy)^*\\vx&\vt_v(a)}
}\eeqn
defines a homomorphism between the linking algebras $\sMatrix{\cB&E^*\\E&\sB^a(E)}$ and $\sMatrix{\vp(\cB)&(vE)^*\\vE&\sB^a(vE)}$ of $E$ and of $vE$, respectively, so that the corner $vE=\Phi_v(E)$ is complete. By this condition the action of $\vt_v(a)$ on $vE$ is determined uniquely as $\vt_v(a)vx=vax$, If $\vt_v$ is well-defined, then it is obviously multiplicative. And the computation
\beqn{
\AB{vx,\vt_v(a^*)vy}
~=~
\AB{vx,va^*y}
~=~
\vp(\AB{x,a^*y})
~=~
\vp(\AB{ax,y})
~=~
\AB{vax,vy}
~=~
\AB{\vt_v(a)vx,vy}
}\eeqn
shows not only that $\vt_v$ is a \nbd{*}map, if it is well-defined, but also that $\vt_v(a)$ is, indeed, well-defined.
\eob

We see that generalized isometries  correspond to homomorphisms of the linking algebra of $E$. Note that the restrictions to the corners of a homomorphism from the linking algebra onto a \nbd{C^*}algebra decompose that \nbd{C^*}algebra into blocks such that the restriction to the corner $E$ becomes a generalized \nbd{\vp}isometry, where $\vp$ is the restriction to the corner $\cB$. In a sense, the correspondence of surjective generalized isometries with $\vp$ also surjective and surjective homomorphisms is one-to-one when $E$ is separable. In that case the restriction of $\vt_v\colon\sK(E)\rightarrow\sK(vE)$ to the compacts extends to a unique surjective homomorphism of the multiplier algebras $\sB^a(E)\rightarrow\sB^a(vE)$; see Pedersen \cite[Proposition 3.12.10]{Ped79} for the extension of homomorphisms and Kasparov \cite[Theorem 1]{Kas80} for $M(\sK(E))=\sB^a(E)$. A condition like separability of $E$ is necessary. (Indeed, following \cite[Example 3.12.11]{Ped79}, if $\Om$ is a nonnormal locally compact Hausdorff space, then there exists a closed subset $A\subset\Om$ and a function in $\sC_b(A)=M(C_0(A))$ that does not admit a continuous extension to $\Om$. Let us put $\cB=E=\sC_0(\Om)$, $\cC=F=\sC_0(A)$ and $\vp(f)=vf=f\upharpoonright A$. Then $v$ is surjective but the homomorphism $\Phi_v$ associated with the surjective \nbd{\vp}isometry is not surjective.) The problem disappears if we pass to von Neumann or \nbd{W^*}modules.

\bprop
Let $\vp$ be a surjective normal homomorphism of von Neumann algebras or \nbd{W^*}algebras and let $v$ be a \nbd{\vp}isometry between von Neumann modules or \nbd{W^*}modules. Then $\Phi_v$ is surjective.
\eprop

\proof
We treat only the case where $\cC$ is a von Neumann algebra, that is, $\cC$ acts nongedenerately as an algebra of operators on a Hilbert space $K$. Then $\sB^a(vE)$ is identified as a concrete algebra of operators acting on the Hilbert space $L:=vE\odot K$. The homomorphism $\sB^a(E)\xrightarrow{\vt_v}\sB^a(vE)\rightarrow\sB(L)$ is normal. (Indeed, $\AB{vx,\vt_v(a)vx}=\AB{vx,vax}=\vp(\AB{x,ax})$ what shows order continuity. See Bhat and Skeide \cite[Appendix C]{BhSk00} or Skeide \cite[Section 3.3]{Ske01} for some explicit calculations of that type.) So the range of this homomorphism is a von Neumann subalgebra of $\sB(L)$. This subalgebra contains $\sK(vE)$ and is generated by $\sK(vE)$ in the strict topology, which is stronger than the strong topology. So the image of $\sB^a(E)$ is at least as big as $\sB^a(vE)$, so, as $\vt_v$ maps into $\sB^a(vE)$ we have $\vt_v(\sB^a(E))=\sB^a(vE)$.\qed

\brem
The preceding proof is minimal in order to obtain the desired result. Actually, one may show that $\Phi_v$ is normal and, therefore, $vE$ is also a strongly closed subset of $\sB(K,L)$. In other words, $vE$ is a von Neumann module.
\erem

\bob
By comparing $\AB{x,a(y+zb)}$ with $\AB{x,ay}+\AB{x,az}\vp(b)$ with the help of the defining equation $\AB{ax,y}=\vp(\AB{x,a^*y})$, we find that a \nbd{\vp}ad\-join\-table map is \nbd{\vp}linear. Also, if $a'$ is another \nbd{\vp}adjointable map, then $a'a^*$ is a usual adjointable map with adjoint $a{a'}^*$ (because $\AB{x,a'a^*y}=\vp(\AB{{a'}^*x,a^*y})=\AB{a{a'}^*x,y}$). Moreover, if $\vp$ is injective on $\cB_E$, then we may apply a left inverse of $\vp$ to the defining equation. We conclude that in this case the \nbd{\vp}adjoint $a^*$ is unique and that also $a^*a'$ is a usual adjointable map with adjoint ${a'}^*a$. Further, $a^*$ is closeable and, therefore, bounded by the \it{closed graph theorem}. If we restrict to the Hilbert $\vp(\cB)$--sub\-mod\-ule $\ol{aE}$ of $F$ (or, otherwise, assume that $\vp$ is also surjective), then $a^*$ is, clearly, \nbd{\vp^{-1}}linear and $a$ its \nbd{\vp^{-1}}adjoint.
\eob

\bprop
Suppose $\vp$ is injective. If $E$ is self-dual (for instance, if $E$ is a von Neumann or \nbd{W^*}module), then every $a\in\sB^{\vp}(E,F)$ has an adjoint.
\eprop

\proof
For simplicity assume that $\vp$ is also surjective (otherwise restrict $\cC$ to $\vp(\cB)\subset\cC$). Recall that $E$ being \hl{self-dual} means that for every bounded right linear map $R\colon E\rightarrow\cB$ there exists a (unique) $x\in E$ such that $R(y)=\AB{x,y}$ for all $y\in E$. (\nbd{W^*}Modules are self-dual by definition, while for von Neumann modules this has been shown in \cite{Ske00b,Ske05c}.) Suppose that $a$ is bounded and \nbd{\vp}linear. Then for every $y\in F$ the map $x\mapsto\vp^{-1}(\AB{y,ax})$ from $E$ to $\cB$ is bounded and right linear. So, there is a unique element $a^*y\in E$ such that $\vp(\AB{a^*y,x})=\vp\bfam{\vp^{-1}(\AB{y,ax})}=\AB{y,ax}$. In oher words, $a^*\colon y\mapsto a^*y$ is a \nbd{\vp}adjoint of $a$.\qed

\bob
Every \nbd{\vp}unitary (ore more generally, every invertible \nbd{\vp}isometry) $u$ has $u^{-1}$ as adjoint. However, since not even usual isometries (that is \nbd{\id}isometries) need to be adjointable, we see that not all \nbd{\vp}linear maps possess a \nbd{\vp}adjoint. More precisely, as for usual isometries one shows that an adjointable \nbd{\vp}isometry $v$ necessarily has complemented range $vE=vv^*F$ in $F$. (Just take into account  that by $\AB{x,vv^*vv^*y}=\AB{vv^*x,vv^*y}=\vp(\AB{v^*x,v^*y})=\AB{x,vv^*y}$ the self-adjoint operator $vv^*$ is an idempotent, so that $(\id_F-vv^*)F$ is the complement of $vE$.) If $\vp\upharpoonright\cB_E$ is injective (and, therefore, also $v$), then the converse statement is also true. (Use the facts that $v$ has a \nbd{\vp}adjoint when considered as mapping onto $pF$, and that the canonical injection $pF\rightarrow F$ has an adjoint. See, for instance, \cite[Proposition 1.5.13]{Ske01} for details.) In general, not even the canonical injection $i_\vp$ need be \nbd{\vp}adjointable. For instance, if $\vp$ is the canonical injection of a subalgebra $\cB$ of $\cC$ into $\cC$, then $i_\vp=\vp$ is \nbd{\vp}adjointable, if and only if $\cB=\cC$. To see necessity assume $\U\in\cB$ (otherwise, use approximate units) and calculate $c=\U^*c=\AB{\vp(\U),c}=\vp(\AB{\U,\vp^*(c)})\in\vp(\cB)=\cB$.
\eob

\bob\label{funcob}
Clearly, the composition $a=a_1a_2$ of two \nbd{\vp_i}linear maps $a_i$ is a \nbd{\vp_1\circ\vp_2}linear map. The same observation holds for generalized isometries and for generalized unitaries. Therefore, under conditions where $a$ determines $\vp$ uniquely, we have a sort of grading on the corresponding sets of generalized maps. This is, in particular, the case, when we restrict our attention to generalized unitaries and \hl{full} modules, that is, to modules $E$ for which $\cB_E=\cB$.
\eob

For every nondegenerate homomorphism $\vp\colon\cB\rightarrow\cC$ and every Hilbert \nbd{\cB}module $E$ there exists a \nbd{\vp}isometry, namely, $i_\vp$. It is injective and surjective, if $\vp\upharpoonright\cB_E$ is injective and surjective, respectively, where for injectivity the condition is also necessary.

We show that an arbitrary \nbd{\vp}linear map factors into the canonical map $i_\vp$ and a usual module map, and we draw some consequences.

\bprop
Let $\vp\colon\cB\rightarrow\cC$ be a nondegenerate homomorphism and let $E$ and $F$ be a Hilbert \nbd{\cB}module and a Hilbert \nbd{\cC}module, respectively, and  denote by $i_\vp$ the canonical map $E\rightarrow E\odot{_\vp\cC}$. Choose $a\in\sB^\vp(E,F)$. Then
\beqn{
a'
\colon
x\odot c
~\longmapsto~
(ax)c
}\eeqn
defines a \nbd{\cC}linear operator from the algebraic tensor product $E\,\ul{\odot}~{_\vp\cC}$ into $F$ such that $a=a'i_\vp$ (uniquely determined by this property if $i_\vp$ is surjective). $a'$ extends to an element in $\sB^r(E\odot{_\vp\cC},F)$ (also denoted by $a'$) with norm $\norm{a'}\le\norm{a}$, if $a$ is \nbd{\vp}adjointable (in that case $a'$ is adjointable) or if $\vp(\cB)$ is an ideal in $\cC$.
\eprop

\proof
Suppose $a'$ is well-defined. If $\cB$ is unital (and, therefore, $\vp(\U_\cB)=\U_\cC$), then $a'i_\vp x=a'(x\odot\U_\cC)=ax\U_\cC=ax$. If $\cB$ is nonunital, then we may use an approximate unit for $\cB$ and, taking into account also Observation \ref{auob}, a similar computation yields again $a'i_\vp x=ax$. Of course, $a'$ is also \nbd{\cC}linear.

The mapping $(x,c)\mapsto(ax)c$ is balanced over $\cB$, that is, $(xb,c)$ and $(x,bc)$ are mapped to the same element. (To verify this, compute $(a(xb))c=(ax)\vp(b)c=(ax)(bc)$.) This alone shows that $a'$ is well-defined on the algebraic tensor product $E\,\ul{\odot}~{_\vp\cC}$. (Anyway, if we show that $a'$ has a (formal) adjoint or that it is bounded, also this will prove well-definedness.)

Next, we compute
\beqn{
\textstyle
\Bnorm{a'\Bfam{\sum_ix_i\odot\vp(b_i)}}
=
\Bnorm{\sum_i(ax_i)\vp(b_i)}
=
\Bnorm{\sum_ia(x_ib_i)}
\le
\norm{a}\Bnorm{\sum_ix_ib_i}
=
\norm{a}\Bnorm{\sum_ix_i\odot\vp(b_i)},
}\eeqn
which shows that the restriction of $a'$ to $E\odot\vp(\cB)$ is bounded by $\norm{a}$. So, if $\vp$ is surjective, then we are done. More generally, if we can show that the range of this restriction (contained in $\cls F\vp(\cB)$) is a pre-Hilbert \nbd{\vp(\cB)}submodule of $F$, that is, if we can show $\AB{ax,ay}\in\vp(\cB)$ for all $x,y\in E$, then $a'=(a'\upharpoonright E\odot\vp(\cB))\odot\id_\cC$ is bounded by $\norm{a}$ as amplification of a mapping in $\sB^r(E\odot\vp(\cB),\cls F\vp(\cB))$ to a mapping $E\odot{_\vp\cC}=(E\odot\vp(\cB))\odot\cC\rightarrow F\odot\cC=F$. Clearly, $\AB{ax,ay}\in\vp(\cB)$, if $\vp(\cB)$ is an ideal in $\cC$ or if $a$ has a \nbd{\vp}adjoint. It remains to note that if $\AB{x,ay}=\vp(\AB{a^*x,y}$, then ${a'}^*\colon xc\mapsto a^*x\odot c$ defines an adjoint of $a'$, because $\AB{y\odot c,a^*x\odot c'}=c^*\AB{y,a^*x}c'=c^*\AB{ay,x}c'=\AB{a'(y\odot c),xc'}$.\qed

\bob
We see that if $i_\vp$ is surjective, then there is a one-to-one correspondence between elements of $\sB^{\vp,a}(E,F)$ and elements in $\sB^a(E\odot{_\vp\cC},F)$. If $\vp(\cB)$ is an ideal, then this correspondence extends to $\sB^\vp(E,F)$ and $\sB^r(E\odot{_\vp\cC},F)$, respectively. All these assertions have much simpler proofs when restricted to \nbd{\vp}unitaries, so that $i_\vp$ and $\vp$ are bijections. In fact, $i_{\vp^{-1}}i_\vp=\id_E$ under the canonical isomorphism $(E\odot{_\vp\cC})\odot{_{\vp^{-1}}\cB}=E\odot({_\vp\cC}\odot{_{\vp^{-1}}\cB})=E\odot\cB=E$. In particular, in this weak sense $i_{\vp^{-1}}$ is the \nbd{\vp}adjoint of $i_\vp$. Suppose that $\vp=\vp_1\circ\vp_2$ for isomorphisms $\vp_i$. Then for every \nbd{\vp}unitary $u$ and \nbd{\vp_i}unitaries $u_i$, the bijections $v_2:=u_1^*u$ and $v_1:=uu_2^*$ are the unique \nbd{\vp_2}unitary and \nbd{\vp_1}unitary, respectively, such that $u_1v_2=u$ and $v_1u_2=u$. In particular, if $u\colon E\rightarrow F$ is a \nbd{\vp}unitary, then $v_E:=ui_{\vp}^*\in\sB^a(E\odot{_\vp}\cC,F)$ and $v_F:=i_{\vp^{-1}}^*u\in\sB^a(F\odot{_{\vp^{-1}}\cB},E)$ are the unique unitaries such that $v_Ei_{\vp}=u=i_{\vp^{-1}}v_F$.
\eob

\bcor\label{phiuucor}
If $u\colon E\rightarrow F$ is a \nbd{\vp}unitary, then $E\odot{_\vp\cC}$ and $F$ are isomorphic. Moreover, $v\longleftrightarrow v_u:=vu^*$ establishes a one-to-one correspondence between \nbd{\vp}unitaries $v\colon E\rightarrow F$ and unitaries $v_u\in\sB^a(F)$, in such a way that $v_uu=v$.
\ecor

In the sequel, we will be particularly interested in the case $F=E$. In this case, Corollary \ref{phiuucor} asserts two things: Firstly, if there exists a \nbd{\vp}unitary in $\sB^{\vp,a}(E)$, then necessarily $E\odot{_\vp\cB}$ is isomorphic to $E$. Secondly, if there exists a \nbd{\vp}unitary in $\sB^{\vp,a}(E)$, then every other \nbd{\vp}unitary in $\sB^{\vp,a}(E)$ factors into the given one and a unique unitary in $\sB^a(E)$. Given a Hilbert \nbd{\cB}module $E$, for an automorphism $\vp$ of $\cB$ we denote by $\SB{\vp}_E$ the set of all automorphisms $\vp'$ that coincide with $\vp$ on $\cB_E$ and we define $\Phi_E:=\bCB{\,\SB{\vp}_E\colon\text{there exists a \nbd{\vp}unitary on $E$}}$.

\bprop\label{phiuniprop}
Suppose $u$ is a \nbd{\vp}unitary on $E$ and let $\vp'$ be another automorphism of $\cB$. Then $u$ is also a \nbd{\vp'}unitary, if and only and only if $\SB{\vp}_E=\SB{\vp'}_E$.
\eprop

\proof
Since the elements $\AB{x,y}$ generate $\cB_E$, the assertion is fairly obvious from $\AB{ux,uy}=\vp^{(')}(\AB{x,y})$, if $u$ is \nbd{\vp^{(')}}unitary.\qed

\bcor
$\Phi_E$ is a group under the product $\SB{\vp}_E\SB{\vp'}_E=\SB{\vp\circ\vp'}_E$ with identity $\SB{\id_\cB}_E$ and inverse $\SB{\vp}_E^{-1}=\SB{\vp^{-1}}_E$.
\ecor

\proof
This follows directly from Observation \ref{funcob}. Indeed, let $\vp_1$ and $\vp_2$ be automorphisms of $\cB$ that admit \nbd{\vp_i}unitaries $u_i$ on $E$, and let $\vp'_i$ be other representatives of $\SB{\vp_i}_E$. Then by Observation \ref{funcob} $u_1u_2$ is both a \nbd{\vp_1\circ\vp_2}unitary and a \nbd{\vp'_1\circ\vp'_2}unitary. By Proposition \ref{phiuniprop} this implies that $\SB{\vp_1\circ\vp_2}_E=\SB{\vp'_1\circ\vp'_2}_E$ so that the product in $\Phi_E$ is well-defined. The rest is obvious.\qed

\brem
Corollary \ref{BEinvcor} will provide a different method to show \it{en passant} that $\Phi_E$ is a group.
\erem

\bcor
$U^{gen}(E)/U(E)=\Phi_E$.
\ecor

We close this section with some considerations that are related to the circle of problems which motivate these notes.

By definition for every $\SB{\vp}_E\in\Phi_E$ there is a \nbd{\vp}unitary $u_\vp$. A natural question is whether we may choose $u_\vp$ in such a way that $u_{\vp_1}u_{\vp_2}=u_{\vp_1\circ\vp_2}$ or, in other words, whether there exists a (necessarily injective) homomorphism $\gamma\colon\Phi_E\rightarrow U^{gen}(E)$ such that $\Phi_E\xrightarrow{\gamma}U^{gen}(E)\rightarrow U^{gen}(E)/U(E)$ is the identity of $\Phi_E$. If the answer is affirmative, then $U^{gen}(E)$ is a semidirect product of $U(E)$ and $\Phi_E$. Indeed, every \nbd{\vp}unitary can be written as $vu_\vp$ with a unique unitary $v$, and if $v_1u_{\vp_1}$ and $v_2u_{\vp_2}$ are an arbitrary \nbd{\vp_1}unitary and \nbd{\vp_2}unitary, respectively, then their product is
\beqn{
v_1u_{\vp_1}v_2u_{\vp_2}
~=~
v_1(u_{\vp_1}v_2u_{\vp_1}^*)u_{\vp_1}u_{\vp_2}
~=~
v_1\alpha_{\vp_1}(v_2)u_{\vp_1\circ\vp_2},
}\eeqn
where $\SB{\vp}_E\mapsto\alpha_\vp:=u_\vp\bullet u_\vp^*$ is a homomorphism $\Phi_E\rightarrow aut(U(E))$. More generally, if $G$ is a subgroup of $\Phi_E$ and $\gamma\colon G\rightarrow U^{gen}(E)$ a homomorhism such that $G\xrightarrow{\gamma}U^{gen}(E)\rightarrow U^{gen}(E)/U(E)$ is the identity of $G$, then the subgroup of $U^{gen}(E)$ generated by $\gamma(G)$ and $U(E)$ is (isomorphic to) a semidirect product of $U(E)$ and $G$.  Moreover, if $\gamma'$ is another such homomorphism, then the two of them are related by a \hl{unitary left cocycle} in the follwoing way. For everey $\vp$ there is a unique unitary $v(\vp)\in U(E)$ such that $\gamma'(\vp)=v(\vp)\gamma(\vp)$. Then $\gamma'$ is a homomorphism, if and only if
\bmun{
\gamma'(\vp_1)\gamma'(\vp_2)
~=~
v(\vp_1)\gamma(\vp_1)v(\vp_2)\gamma(\vp_2)
~=~
v(\vp_1)\alpha_{\vp_1}(v(\vp_2))\gamma(\vp_1\circ\vp_2)
\\
~=~
v(\vp_1\circ\vp_2)\gamma(\vp_1\circ\vp_2)
~=~
\gamma'(\vp_1\circ\vp_2),
}\emun
that is, if and only if $v(\vp_1)\alpha_{\vp_1}(v(\vp_2))=v(\vp_1\circ\vp_2)$. In most examples where $\Phi_E$ can be computed easily, this is so, because it is easy to find a \it{canonical} \nbd{\vp}unitary to every candidate $\vp$ and the map that sends $\SB{\vp}_E$ to that canonical \nbd{\vp}unitary turns out to be a homomorphism. (Among these examples, there are the cases where $E=\cB$ and the cases where $\cB$ is commutative or finite-dimensional.) In general, we do not know the answer.

Among the subgroups of $\Phi_E$ the one-parameter groups $G=\bfam{\vp_t}_{t\in\R}$ are particularly interesting. Excluding the periodic case, so that $t\mapsto\vp_t$ is injective, the possible homomorphisms $\gamma\colon G\rightarrow U^{gen}(E)$ are exactly those one-parameter groups $\bfam{u_t}_{t\in\R}$ in $U^{gen}(E)$ with $u_t\in\sB^{a,\vp_t}(E)$. It follows immediately that two such one-parameter groups differ by a cocycle in $U(E)$. Of course, this remains true also in the general case, that is, if $G$ is not necessarily nonperiodic. Abbaspour, Moslehian and Niknam \cite{AMN05p1} showed that the generators of strongly continuous one-parameter groups of generalized unitaries are generalized \nbd{E}valued derivations. The study of these groups is one of our motivations. We will return to this setting in the end of Section \ref{combSEC}.

The fact that existence of a \nbd{\vp}unitary factors $E$ into $E\odot{_\vp\cB}$ (up to unitary isomorphism) reminds us of the fact that by the representation theory of $\sB^a(E)$ a strict unital endomorphism $\vt$ of $\sB^a(E)$ factors $E$ into $E\odot E_\vt$ where $E_\vt$ is a correspondence over $\cB$; see Muhly, Skeide and Solel \cite{MSS06}. We discuss this in the following section as a preparation for Section \ref{combSEC}, where we try to use the automorphisms of $\sB^a(E)$ induced by a \nbd{\vp}unitary on $E$ to understand better the structure of $\Phi_E$.

\section{Automorphisms and their relation with the Picard group}\label{PicSEC}

{
We start by recalling briefly what \cite[Theorem 1.4]{MSS06} asserts about strict representations of $\sB^a(E)$. Let $\cB$ and $\cC$ be \nbd{C^*}algebras, let $E$ be a Hilbert \nbd{\cB}module and  suppose $\vt\colon\sB^a(E)\rightarrow\sB^a(F)$ is a unital strict representation of $\sB^a(E)$ on a Hilbert \nbd{\cC}module $F$. Here, (among other equivalent descriptions) $\vt$ being unital and strict means that the action of $\sK(E)$ on $F$ (via $\vt$) is nondegenerate. Therefore, $F$ is not only a correspondence from $\sB^a(E)$ to $\cC$ (which we denote by $_\vt F$) but even a correspondence from $\sK(E)$ to $\cC$. Further, $E^*$ with inner product $\AB{x^*,y^*}:=xy^*$ is a correspondence from $\cB$ to $\sB^a(E)$ (with module operations $bx^*a:=(a^*xb^*)^*$) which may also be considered as a correspondence from $\cB_E$ to $\sK(E)$ (actually, a Morita equivalence; see below). It is clear that
\beqn{
F
~\cong~
\sK(E)\odot{_\vt F}
~\cong~
(E\odot E^*)\odot{_\vt F}
~\cong~
E\odot(E^*\odot{_\vt F}).
}\eeqn
The following theorem fixes an isomorphism and summarizes some more results from \cite{MSS06}.

\bitemp[Theorem \cite{MSS06}.~]\label{MSSthm}
Define the correspondence $F_\vt:=E^*\odot{_\vt F}$ from $\cB$ to $\cC$. Then $F\cong E\odot F_\vt$ and $\vt(a)=a\odot\id_{F_\vt}$. More precisely, $u\colon x\odot(y^*\odot z)\mapsto\vt(xy^*)z$ $(x,y\in E;z\in F)$ defines a unitary $E\odot(E^*\odot{_\vt F})\rightarrow F$ such that $\vt(a)=u(a\odot\id_{F_\vt})u^*$.

Moreover, if $F'$ is another correspondence from $\cB$ to $\cC$ which is also a correspondence from $\cB_E$ to $\cC$ (that is, $\sB_E$ acts nondegenerately on $F'$) with a unitary $u'\colon E\odot F'\rightarrow F$ fulfilling $\vt(a)=u'(a\odot\id_{F'}){u'}^*$, then there is a unique \hl{isomorphism} of correspondences (that is, a bilinear unitary) $v\colon F_\vt\rightarrow F'$ fulfilling $u'(\id_E\odot v)=u$. In particular, if $E$ is full, then the \hl{multiplicity correspondence} $F_\vt$ is unique up to isomorphism.
\eitemp

A \hl{Morita equivalence} from $\cB$ to $\cC$ is a correspondence $M$ from $\cB$ to $\cC$ such that $\cC_M=\cC$ (that is, $M$ is full) and $\sK(M)=\cB$ (meaning that the canonical action of $\cB$ on $M$ defines an isomorphism onto $\sK(M)$). By these properties it follows that $M^*$ is a Morita equivalence from $\cC$ to $\cB$. We summarize \cite[Corollary 1.11 and Remark 1.12]{MSS06}: If $\vt\colon\sB^a(E)\rightarrow\sB^a(F)$ is an isomorphism that is \hl{bistrict} (that is, both $\vt$ and $\vt^{-1}$ are strict), then $F_\vt$ is a \hl{Morita equivalence} from $\cB_E$ to $\cC_F$ and $F_{\vt^{-1}}=F_\vt^*$. (Anoussis and Todorov \cite[Corollary 2.5]{AnTo05} recently proved that if the involved unital \nbd{C^*}algebras are unital and separable and if the modules are separable, then an isomorphism $\vt$ is bistrict, automatically.)

Morita equivalences compose under tensor product. More precisely, if $M$ is a Morita equivalence from $\cB$ to $\cC$ and $N$ is a Morita equivalence from $\cC$ to $\cD$, then $M\odot N$ is a Morita equivalence from $\cB$ to $\cD$. In particular, $M\odot M^*\cong\cB$ and $M^*\odot M=\cC$, that is, $M$ and $M^*$ are inverses up to isomorphism under this composition. Applied to the composition of bistrict isomorphisms $\sB^a(E)\xrightarrow{\vt}\sB^a(F)\xrightarrow{\theta}\sB^a(G)$, we find $F_\vt\odot G_\theta\cong G_{\theta\circ\vt}$. But, \cite[Theorem 1.14]{MSS06} tells more: Together with the concrete constructions in Theorem \ref{MSSthm} the isomorphisms $F_\vt\odot G_\theta\cong G_{\theta\circ\vt}$ can be chosen such that they iterate associatively.

The isomorphism classes $\SB{M}_\cB$ of Morita equivalences $M$ from $\cB$ to $\cB$ (or over $\cB$) form a group under tensor product, the \hl{Picard group} $Pic(\cB)$ of $\cB$; see Brown, Green and Rieffel \cite{BGR77}. Clearly, if $\vp$ is an automorphism of $\cB$, then $_\vp\cB$ is a Morita equivalence over $\cB$, and $\SB{_\vp\cB}_\cB=\SB{_{\vp'}\cB}_\cB$ in $Pic(\cB)$, if and only if $\vp$ and $\vp'$ differ by a \hl{quasi inner automorphism} of $\cB$, that is by an automorphism which is given by conjugation with a unitary in the multiplier algebra $M(\cB)=\sB^a(\cB)$ of $\cB$. We denote by $gin(\cB)$ the group of quasi inner automorphisms. Since every automorphism of $\cB$ extends to a unique automorphism of $M(\cB)$, the subgroup $gin(\cB)$ of $aut(\cB)$ is normal and by sending $\SB{\vp}_{gin(\cB)}$ to $\SB{_\vp\cB}_\cB$ we define an injective homomorphism $aut(\cB)/gin(\cB)\rightarrow Pic(\cB)^{op}$ (\cite[Proposition 3.1]{BGR77}) where $G^{op}$ denotes the \hl{opposite group} of the group $G$ with multiplication $g*g'=g'g$. By \cite[Corollary 3.5]{BGR77} this homomorphism is surjective, for instance, if $\cB$ is stable and has a strictly positive element. But it need not be surjective. A concrete counter example is the Morita equivalence $M=\tMatrix{0&\C&\C\\\C&0&0\\\C&0&0}$ over $\cB=\tMatrix{\C&0&0\\0&\C&\C\\0&\C&\C}$. Then $M$ does not have a vector $x$ with $\AB{x,x}=\U$ and, therefore, cannot be isomorphic to the unital algebra $\cB$, not even as a right module.

Now we turn our interest to bistrict automorphisms $\vt$ of $\sB^a(E)$ for a fixed not necessarilly full Hilbert \nbd{\cB}module $E$. We denote the group of all these automorphisms by $straut(\sB^a(E))$ and by $inn(\sB^a(E))$ the (clearly, normal) subgroup of inner automorphisms. By Theorem \ref{MSSthm} and the forthcoming discussion we associate with every $\vt\in straut(\sB^a(E))$ a Morita equivalence $E_\vt=E^*\odot{_\vt E}$ over $\sB_E$.

\bprop
Let $\vt_1$ and $\vt_2$ be bistrict automorphisms of $\sB^a(E)$. Then $E_{\vt_1}\cong E_{\vt_2}$, if and only if $\vt_1$ and $\vt_2$ differ by an inner automorphism of $\sB^a(E)$.
\eprop

\proof
We prove only the direction that is not fairly obvious. So, suppose that $u\in\sB^a(E)$ is a unitary such that $\vt_2(a)=u\vt_1(a)u^*$. Let us indicate by $x^*\odot_i y$ an elementary tensor in $E^*\odot{_{\vt_i}E}$. One easily checks that $x^*\odot_1y\mapsto x^*\odot_2uy$ defines an isometry $E_{\vt_1}\rightarrow E_{\vt_2}$ that, clearly, is also surjective and left linear.\qed

\bcor\label{strautcor}
By sending $\SB{\vt}_{inn(\sB^a(E))}$ to $\SB{E_\vt}_{\cB_E}$ we define an injective homomorphism
\beqn{
straut(\sB^a(E))/inn(\sB^a(E))
~\longrightarrow~
Pic(\cB_E)^{op}.
}\eeqn
\ecor

\brem
It would be tempting to consider directly the injective mapping $\vt\rightarrow E_\vt$ without dividing out equivalence classes. But there is no possibility to discuss away the fact that $E_{\vt_2}\odot E_{\vt_1}$ and $E_{\vt_1\circ\vt_2}$ are isomorphic but not equal in the category of correspondences over $\cB_E$. In fact, the Morita equivalences over $\cB_E$ do not form a group under tensor product (not even a monoid!) but only a semigroup. Taking the quotient over $gin(\cB_E)$ gives the group property. But now $\vt\rightarrow\SB{E_\vt}_{\cB_E}$ is no longer injective. In fact, Corollary \ref{strautcor} identifies correctly its kernel as $inn(\sB^a(E))$.
\erem

A Morita equivalence $M$ over $\cB_E$  is isomorphic to $E_\vt$ for some $\vt\in straut(\sB^a(E))$, if and only if there exists a unitary $E\rightarrow E\odot M$. If we put $E=M=\tMatrix{0&\C&\C\\\C&0&0\\\C&0&0}$ (cf.\ the preceding counter example), then $E\odot M=\cB\ncong E$. This shows that also the homomorphism in Corollary \ref{strautcor} need not be surjective. At least, in a variant for von Neumann algebras and von Neumann modules (all homomorphisms being normal) the results from Skeide \cite{Ske04p} assert that one may find for every $\SB{M}_\cB$ in the (von Neumann version of the) Picard group a strongly full von Neumann module $E$ with a normal automorphism $\vt$ of $\sB^a(E)$ such that $M$ is the strong closure of $E^*\odot{_\vt E}$. The example considered in \cite{Ske04p} is exactly $M$ and thanks to the fact that $\cB$ is finite-dimensional the constructed automorphism is also strict.

Both homomorphisms $aut(\cB_E)/gin(\cB_E)\rightarrow Pic(\cB_E)^{op}$ and $straut(\sB^a(E))/inn(\sB^a(E))\rightarrow Pic(\cB_E)^{op}$ are injective. So we may ask, whether the image of one is contained in the image of the other (identifying in that way one of the groups $straut(\sB^a(E))/inn(\sB^a(E))$ and $aut(\cB_E)/gin(\cB_E)$ as a subgroup of the other). However, the preceding discussion shows that even if $E$ is supposed full, the image of $straut(\sB^a(E))/inn(\sB^a(E))\rightarrow Pic(\cB_E)^{op}$ need not be contained in the image of $aut(\cB)/gin(\cB)\rightarrow Pic(\cB)^{op}$. We close with a counter example that shows that also the converse statement need not be true.

\bex\label{phnotthex}
Let $\cB=\cA^2$ for a \nbd{C^*}algebra $\cA$ (for instance, $\cA=\C$ so that $\cB=\C^2$). Define the flip $\vp\sMatrix{a_1\\a_2}=\sMatrix{a_2\\a_1}$. Set $E=\sMatrix{\cA\\\zero}\oplus\sMatrix{\cA\\\zero}\oplus\sMatrix{\zero\\\cA}$. Then $E$ is full and $E\odot{_\vp\cB}=\sMatrix{\zero\\\cA}\oplus\sMatrix{\zero\\\cA}\oplus\sMatrix{\cA\\\zero}$ is not isomorphic to $E$. This shows that $_\vp\cB$ is not isomorphic to $E_\vt$ for any automorphism $\vt$ of $\sB^a(E)$.
\eex
}

\section{Relating $\Phi_E$ and $Pic(\cB_E)$}\label{combSEC}

In this section we apply our knowledge from Section \ref{PicSEC} to the bistrict automorphisms $\vt_u=u\bullet u^*$ induced by conjugation with a \nbd{\vp}unitary $u$ to understand better the group $\Phi_E$ defined in Section \ref{genmodmapSEC}.

Throughout, $E$ is a (not necessarily full) Hilbert \nbd{\cB}module. The first thing to do is to convince ourselves that, indeed, $\vt_u\in straut(\sB^a(E))$ for every \nbd{\vp}unitary $u\in U^{gen}(E)$. As $\vt_u$ is, clearly, multiplicative, it is sufficient to show that $\vt_u$ is a \nbd{*}map from which it follows, too, that $\vt_u$ maps into (and, therefore, onto) $\sB^a(E)$. We compute
\beqn{
\AB{x,\vt_u(a^*)y}
~=~
\AB{x,ua^*u^*y}
~=~
\vp(\AB{u^*x,a^*u^*y})
~=~
\vp(\AB{au^*x,u^*y})
~=~
\AB{uau^*x,y}
~=~
\AB{\vt_u(a)x,y}.
}\eeqn
We observe also that $\vt_u(xy^*)=(ux)(uy)^*$, because
\beqn{
\vt_u(xy^*)z
~=~
u(xy^*)u^*z
~=~
u(x\AB{y,u^*z})
~=~
(ux)\vp(\AB{y,u^*z})
~=~
(ux)\AB{uy,z}
~=~
(ux)(uy)^*z.
}\eeqn
So, $\vt_u$ maps $\sK(E)$ into (and, therefore, onto) $\sK(E)$. In other words, $\vt_u$ is bistrict.

Identifying correctly the correspondence over $\cB_E$ of $\vt_u$ we see that mere existence of a \nbd{\vp}unitary implies that $\vp$ leaves $\cB_E$ (globally) invariant (Corollary \ref{BEinvcor} below).

\bprop\label{geninnprop}
$E_{\vt_u}\cong{_\vp}(\cB_E)$.
\eprop

\proof
The mapping $x^*\odot y\mapsto\AB{ux,y}$ defines an isometry $E^*\odot{_{\vt_u}E}\rightarrow\cB_E$. Indeed, we find
\beqn{
\AB{{x'}^*\odot y',x^*\odot y}
~=~
\AB{y',\vt_u(x'x^*)y}
~=~
\AB{y',(ux')(ux)^*y}
~=~
\AB{y',ux'}\AB{ux,y}.
}\eeqn
Cleary, this isometry is surjective, that is, a unitary. By
\beqn{
b(x^*\odot y)
~=~
(xb^*)^*\odot y
~\longmapsto~
\AB{u(xb^*),y}
~=~
\AB{(ux)\vp(b^*),y}
~=~
\vp(b)\AB{ux,y}
}\eeqn
we see that it is also left linear, that is, an isomorphism of correspondences.\qed

\lf\noindent
This computation is for correspondences over $\cB$. But we know that $E_{\vt_u}$ and, therefore, also ${_\vp}(\cB_E)$ may be viewed as a correspondence over $\cB_E$ and that as such it must be a Morita equivalence over $\cB_E$. In particular, $\vp$ must map $\cB_E$ onto $\sK(\cB_E)=\cB_E$.

\bcor\label{BEinvcor}
$\vp$ (co-)restricts to an automorphism of $\cB_E$.
\ecor

Of course, this follows also directly from the definition of \nbd{\vp}unitary and the fact that the inverse of a \nbd{\vp}unitary is a \nbd{\vp^{-1}}unitary in Section \ref{genmodmapSEC}.

\bex
If, in Example \ref{phnotthex}, we choose $E=\sMatrix{\cA\\\zero}$ then $\cB_E=\sMatrix{\cA\\\zero}$ and $\vp$ does not leave $\cB_E$ invariant. Therefore, there is no \nbd{\vp}unitary on $E$.
\eex

We find that $\Phi_E$ is a subgroup of $aut(\cB_E)$. Namely, in order to determine an element of $\Phi_E$ an automorphism $\vp_E$ of $\cB_E$ must fulfill two conditions. Firstly, $\vp_E$ must admit an extension to an automorphism of $\cB$. (Which one is not important as the class $\SB{\bullet}_E$ ignores differences outside $\cB_E$.) Secondly, $\vp_E$ must admit a \nbd{\vp_E}unitary on $E$. (Since $\cB_E$ is an ideal in $\cB$ and $\vp$ leaves $\cB_E$ invariant, we find $E\odot{_{\vp}\cB}=E\odot{_{\vp_E}\cB_E}$. In other words, if $u$ is a \nbd{\vp_E}unitary and $\vp_E$ admits an extension $\vp$ to $\cB$, then $u$ is also a \nbd{\vp}unitary.)

We see that the problem of deciding whether an automorphism $\vp_E$ of $\cB_E$ is an element of $\Phi_E$ decomposes into two independent questions, which both must have affirmative answers. The first condition arises by considering the full Hilbert \nbd{\cB_E}module $E$ as a Hilbert \nbd{\cB}module for a \nbd{C^*}algebra $\cB$ that contains $\cB_E$ as an ideal. It is an extendability condition which regards only the extension of automorphisms $\vp_E$ of $\cB_E$ to automorphisms $\vp$ of $\cB$. The concrete form of $E$ does not appear except that it determines the ideal $\cB_E$. If $E$ is full the problem disappears. Also if $\cB_E$ is unital (so that $\cB$ decomposes into the direct sum of $\cB_E$ and its complement $\cB_E^\perp$) every automorphism of $\cB_E$ extends (for instance, by $\id_{\cB_E^\perp}$ on the complement) to $\cB$. The same remains true if $\cB$ contains the multiplier algebra of $\cB_E$ as a direct summand (because every automorphism of $\cB_E$ extends to an automorphism of the multiplier algebra). The second condition regards existence of a \nbd{\vp_E}unitary for the full Hilbert \nbd{\cB_E}module $E$. In view of Observation \ref{isolinob} also this is a problem of extendability, now of an automorphism of the corner $\cB_E$ of the linking algebra to an automorphism of the whole linking algebra. We may view the first and the second condition as a single extension problem as follows.

\bprop
An automorphism $\vp_E$ of $\cB_E$ is in $\Phi_E$, if and only if it admits an extension as a matrix \nbd{C^*}algebra automorphism $\Phi$ of the linking algebra $\sMatrix{\cB&E^*\\E&\sB^a(E)}\supset\sMatrix{\cB_E&0\\0&0}\cong\cB_E$. The answer is affirmative if and only if $\vp_E$ admits both an extension as an automorphism to $\cB$ and an extension as a matrix \nbd{C^*}algebra automorphism to $\sMatrix{\cB_E&E^*\\E&\sB^a(E)}$.
\eprop

Recall that according to \cite{Ske00b} a matrix \nbd{C^*} algebra automorphism is a \nbd{C^*}algebra automorphism of a matrix \nbd{C^*}algebra that respects the corners. In Observation \ref{isolinob} this is automatic as the extension of the homomorphism is into a different algebra and defines there a suitable decomposition. Here the range of the extension is given together with a decomposition, so that we must require explicitly that the extension respects the decomposition. Note also that the restriction to $\sB^a(E)$ is automatically of the form $\vt_u$ where $u$ is the restriction of $\Phi$ to $E$ and if $\vp$ is the restriction of $\Phi$ to $\cB$ (that is, an extension of $\vp_E$ to $\cB$) then $u$ is a \nbd{\vp}unitary.

\lf
So far, we have $U^{gen}(E)/U(E)=\Phi_E\subset aut(\cB_E)$ and we have identified $aut(\cB_E)/gin(\cB_E)$ and $straut(\sB^a(E))/inn(\sB^a(E))$ as (possibly different) subgroups of $Pic(\cB_E)^{op}$. We even know that $u\mapsto\vt_u$ defines a homomorphism $U^{gen}(E)\rightarrow straut(\sB^a(E))$ and that $\vt_{U(E)}=inn(\sB^a(E))$. This shows that $\SB{u}_{U(E)}\mapsto\SB{\vt_u}_{inn(\sB^a(E))}$ well-defines a homomorphism $\Phi_E\rightarrow Pic(\cB_E)^{op}$.

\bthm
$\Phi_E\cap gin(\cB_E)$ is a normal subgroup of $\Phi_E$ and the kernel of $\Phi_E\rightarrow Pic(\cB_E)^{op}$ is exactly $\Phi_E\cap gin(\cB_E)$. In particular, we find
\beqn{
\Phi_E/(\Phi_E\cap gin(\cB_E))
~\subset~
straut(\sB^a(E))/inn(\sB^a(E))
~\subset~
Pic(\cB_E)^{op}.
}\eeqn
\ethm

\proof
By Corollary \ref{strautcor} and Proposition \ref{geninnprop} the image of $\vp\in\Phi_E$ in $Pic(\cB_E)$ is the same as the image under the canonical map $aut(\cB_E)\rightarrow Pic(\cB_E)$, namely, $\SB{_\vp\cB_E}_{\cB_E}$. This means that an automorphism $\vp\in\Phi_E$ is in the kernel, if and only if the restriction to $\cB_E$ is in $gin(\cB_E)$, that is, if $\vp$ is in $\Phi_E\cap gin(\cB_E)$. In particular, $\Phi_E\cap gin(\cB_E)$ is the kernel of a homomorphism and, therefore, a normal subgroup. The remaing statements are obvious.\qed

\lf
Once again, in order to decide whether an element of $gin(\cB)$ is in $\Phi_E$ we have to face two problems. Firstly, it is unclear whether an automorphism in $gin(\cB_E)$ extends as an automorphism to $\cB$, and, if it extends, whether this extension can be chosen as an automorphism in $gin(\cB)$. (Even if $\cB_E$ is essential in $\cB$, so that $\cB_E\subset\cB\subset M(\cB)\subset M(\cB_E)$, it is unclear, whether a unitary $u$ in $M(\cB_E)$ fulfills $u\cB u^*=\cB$, and, if it fulfills $u\cB u^*=\cB$, whether $u$ can be chosen in $M(\cB)$.) In general, if $\cB$ contains $M(\cB_E)$ as a direct summand (for instance, if $\cB_E$ is unital), then every element of $gin(\cB_E)$ extends as an element of $gin(\cB)$. Secondly, it is unclear whether for full $E$ an element $\vp_E$ of $gin(\cB_E)$ admits a \nbd{\vp_E}unitary. But, unlike the case of a general automorphism of $\cB_E$, here we are able to show that $\Phi_E\subset aut(\cB_E)$ contains $gin(\cB_E)$ as an (of course, normal) subgroup:

\bprop
Let $v$ be a unitary in $M(\cB_E)$ and set $\vp_v=v\bullet v^*\in gin(\cB_E)$. Choose a bounded approximate unit $\bfam{u_\lambda}_{\lambda\in\Lambda}$ for $\cB_E$. Then $u_v\colon x\mapsto\lim_\lambda x(u_\lambda v^*)$ defines a \nbd{\vp_v}unitary $u_v$.
\eprop

\proof
$xu_\lambda$ converges to $x$ so that it is, in particular, a Cauchy net. In other words, $\bnorm{x((u_\lambda-u_{\lambda'})v^*)}^2\le\bnorm{x(u_\lambda-u_{\lambda'})}^2\le\ve$ for $\lambda,\lambda'$ big enough. It is routine to check that $u_v$ is a \nbd{\vp_v}isometry. To show that it is surjective, we simply observe that $u_{v^*}$ is an inverse.\qed

\lf
We note that $v\mapsto u_v$ is an injective homomorphism.

\bcor
If $E$ is a full Hilbert \nbd{\cB}module, then $U(M(\cB))\subset U^{gen}(E)$. The image in $\Phi_E$ is exactly $gin(\cB)$.
\ecor

\brem
We mentioned already in the end of Section \ref{genmodmapSEC} that our motivation to study the generalized unitary group comes from generalized dynamical systems on Hilbert modules, that is, strongly continuous one-parameter groups of generalized unitaries as studied in \cite{AMN05p1}. We think that our analysis has a large potential for this sort of problemes. We showed already in the end of Section \ref{genmodmapSEC} that two such groups are conjugate by a cocycle in $U(E)$, if and only if their images in $\Phi_E$ under the canonical map $U^{gen}(E)\rightarrow\Phi_E$ coincide. In other words, we discussed perturbations by unitary cocylces of automorphism groups on $\sB^a(E)$ induced by generalized unitary groups. It is natural to apply the results from \cite{AMN05p1} about the generators to this case. But, first, this is too far reaching for this introduction and, second, the results of \cite{AMN05p1} concern only  one side of the coin, namely, they assert that a generator of a generalized unitary group is a generalized derivation. Sufficient conditions on a generalized derivation to be a generator are yet missing.

Of course, if $u_t$ is a \nbd{\vp_t}unitary and if $v_t$ is a \nbd{\psi_t}unitary, then $u_tv_tu_t^*$ is also a \nbd{\psi_t}unitary. Therefore, it has sense to ask for the perurbation of the automorphism group $\bfam{u_t\bullet u_t}_{t\in\R}$ by a cocylce of generalized unitaries $v_t$. One obtains a cocycle condition for the $\psi_t$ with respect to the group $\bfam{\vp_t}_{t\in\R}$ and a cocycle condition for $v_t$ which is completely parallel to the case of unitary cocycles.

Far more involved is the question for perturbations of general strict \nbd{E_0}semigroups $\vt=\bfam{\vt_t}_{t\in\R_+}$ (that is, semigroups of unital endomorphisms $\vt_t$) on $\sB^a(E)$ by families $\bfam{v_t}_{t\in\R_+}$ of generalized unitaries. In other words, we are seeking conditions on the generalized unitaries $v_t$ that guarantee that the maps $\vt_t^v=v_t\vt_t(\bullet)v_t^*$ still form an \nbd{E_0}semigroup. For unitaries $v_t$ one obtains the well-known cocycle condition $v_{s+t}=v_s\vt_s(v_t)$. For the generalization of this condition to generalized unitaries we have to face (at least) two problems. The more obvious one is that we have to give a meaning to $\vt_s(v_t)$ for a generalized unitary $v_t$, because so far $\vt_s$ is defined only for elements of $\sB^a(E)$. The more hidden problem is that in the verification that an equation like $v_{s+t}=v_s\vt_s(v_t)$ suffices for that the $\vt_t^v$ form a semigroup, one has to verify that the extension of $\vt_s$ to elements of the form $v_tav_t^*$ $(a\in\sB^a(E))$ behaves a sort of multiplicatively. For these problems, so far, we did not yet spot an obvious solution. It is here where we believe that our discussion, putting so much emphasis on the role played by the tensor product, will be crucial. \nbd{E_0}Semigroups come along with product systems $E^\odot=\bfam{E_t}_{t\in\R_+}$. It is natural to try to look for a substitute for the cocycle condition in terms of the associated product systems and we are convinced that conjugation of the correpsondences $E_t$ with $_{\vp_t}\cB$ (where $v_t$ is a \nbd{\vp_t}unitary), that is, expressions like $_{\vp_t}\cB\odot E_t\odot{_{\vp_t^{-1}}\cB}$, will play an outstanding role. But this must be rerserved for future work.
\erem

\setlength{\baselineskip}{2.5ex}


\newcommand{\Swap}[2]{#2#1}\newcommand{\Sort}[1]{}
\providecommand{\bysame}{\leavevmode\hbox to3em{\hrulefill}\thinspace}
\providecommand{\MR}{\relax\ifhmode\unskip\space\fi MR }
\providecommand{\MRhref}[2]{%
  \href{http://www.ams.org/mathscinet-getitem?mr=#1}{#2}
}
\providecommand{\href}[2]{#2}


\end{document}